\documentclass[11pt, a4paper]{article}
\usepackage{graphicx}
\usepackage{amsfonts}
\usepackage{amsmath}
\usepackage{amssymb}
\usepackage{amsthm}
\usepackage{bbm}
\usepackage[ansinew]{inputenc}
\usepackage{latexsym}
\usepackage{tabularx}
\allowdisplaybreaks
\newtheorem{thm}{Theorem}
\newtheorem{cor}{Corollary}

\newtheorem{defn}{Definition}
\newtheorem{rem}{Remark}

\newtheorem{examp}{Example}

\newcommand{\norm}[1]{\left\Vert#1\right\Vert}
\newcommand{\abs}[1]{\left\vert#1\right\vert}

\newcommand{\To}{\rightarrow}

\newcommand{\bsgamma}{\boldsymbol{\gamma}}
\newcommand{\bsdelta}{\boldsymbol{\delta}}

\newcommand{\bsbeta}{\boldsymbol{\beta}}
\newcommand{\bseta}{\boldsymbol{\eta}}

\newcommand{\bsx}{\boldsymbol{x}}

\newcommand{\bse}{\boldsymbol{e}}

\newcommand{\bsc}{\boldsymbol{c}}

\newcommand{\bstau}{\boldsymbol{\tau}}
\newcommand{\bsid}{\boldsymbol{\rm id}}

\newcommand{\bsi}{\boldsymbol{i}}

\newcommand{\bssigma}{\boldsymbol{\sigma}}

\newcommand{\cH}{{\cal H}}

\newcommand{\id}{{\rm id}}

\newcommand{\bszero}{\boldsymbol{0}}

\newcommand{\rd}{\,\mathrm{d}}

\newcommand{\NN}{\mathbbm{N}}
\newcommand{\ZZ}{\mathbbm{Z}}

\newcommand{\Sy}{\mathfrak{S}}

\makeatletter
\newcommand{\rdots}{\mathinner{\mkern1mu\lower-1\p@\vbox{\kern7\p@\hbox{.}}
\mkern2mu \raise4\p@\hbox{.}\mkern2mu\raise7\p@\hbox{.}\mkern1mu}}
\makeatother

\date{}

\begin{document}

\title{Discrepancy bounds for low-dimensional point sets}

\author{Henri Faure and Peter Kritzer\thanks{P.~Kritzer gratefully
acknowledges the support of the Austrian Science Fund (FWF), Project
P23389-N18, and Project F5506-N26, which is part of the Special Research Program ``Quasi-Monte
Carlo Methods: Theory and Applications''}}

\maketitle

\begin{center}
Dedicated to H.~Niederreiter on the occasion of his 70th birthday.
\end{center}

\begin{abstract}
\noindent The class of $(t,m,s)$-nets and $(t,s)$-sequences, introduced in their 
most general form by Niederreiter, are important examples of point 
sets and sequences that are commonly used in quasi-Monte Carlo 
algorithms for integration and approximation. Low-dimensional versions 
of $(t,m,s)$-nets and $(t,s)$-sequences, such as Hammersley point sets 
and van der Corput sequences, form important sub-classes, as they are 
interesting mathematical objects from a theoretical point of view, and 
simultaneously serve as examples that make it easier to understand the 
structural properties of $(t,m,s)$-nets and $(t,s)$-sequences in 
arbitrary dimension. For these reasons, a considerable number of 
papers have been written on the properties of low-dimensional nets and 
sequences.

In this paper, we summarize recent results on the distribution 
properties of low-dimensional examples of $(t,m,s)$-nets and 
$(t,s)$-sequences,  state a new result regarding lower star discrepancy bounds, and
formulate some open questions.
\end{abstract}



\section{Introduction}\label{Intro}

In many applications of mathematics, such as in finance or computer graphics, one is in need of numerically approximating the value of an integral 
$I_s(f)=\int_{[0,1]^s} f(\bsx)\rd\bsx$
of a function $f$ defined on $[0,1]^s$. 
One way of dealing with this problem is to use a \textit{quasi-Monte Carlo} integration rule, which is an equal weight 
quadrature rule of the form $Q_{N,s}(f):=N^{-1}\sum_{n=0}^{N-1} f(\bsx_n)$, 
where the integration nodes $\bsx_0,\ldots,\bsx_{N-1}$ are deterministically chosen points in the unit cube. 
It is well known in the theory of quasi-Monte Carlo methods
that a useful property of the points $\bsx_0,\ldots,\bsx_{N-1}$ is that 
they are very evenly distributed in the integration domain $[0,1]^s$. In this context, we frequently refer to  
the points $\bsx_0,\ldots,\bsx_{N-1}$ as a point set, by which we mean a multi-set, i.e., points may occur repeatedly. 
In general, the  term ``point set'' also includes infinite sets, i.e., infinite sequences. 

In order to measure uniformity of distribution of a given point set $P$ with points $\bsx_{0}, \ldots,\bsx_{N-1}$ in $[0,1)^{s}$, one frequently studies the star discrepancy $D^*$, which is defined by
\[D^*(N, P)=\sup_{J}\abs{A(J,N,P)-N \lambda (J)}:=\sup_J \abs{\Delta(J,N,P)},\]
where the supremum is extended over all intervals $J\subseteq[0,1)^{s}$ of the form
$J=\prod_{j=1}^{s}[0,\alpha_{j})$, $0<\alpha_{j}\le 1$, $A(J,N,P)$ denotes the number of $i$ with $\bsx_{i}\in J$,
and $\lambda$ is the Lebesgue measure ($\Delta (J,N,P)$ is the so-called local discrepancy function).

In one dimension, where very precise results exist, we will also need the notion 
of (extreme) discrepancy $D$ obtained by taking the supremum of $\abs{\Delta (J,N,P)}$ 
over all intervals $J$ (not necessarily anchored in the origin).

When considering an infinite sequence $S=(\bsx_n)_{n\ge 0}$ of points in $[0,1)^s$, 
we denote by $D^*(N, S)$ (resp.~$D(N,S)$) the star discrepancy (resp.~the discrepancy) of the first $N$ elements of $S$. 
In the case of a finite point set $P$ with $N$ points, we frequently write $D^* (P)$ (resp.~$D(P)$) if there 
is no possible confusion regarding the number of points.

{\bf Relation between sequences and finite point sets.}
A general principle (also valid in arbitrary dimension) states the link between one-dimensional sequences and two-dimensional point sets deduced from them 
\cite{F86, N92}: let $S=(x_n)_{n\ge 0}$ be an infinite sequence taking its values in $[0,1]$ and let 
$P$ be the two-dimensional point set  
$P=\left \{ \left (x_n, \frac{n}{N} \right )\ ;\ 0 \leq n < N \right \}\subset [0,1]^2. \ \mbox{ Then }$
\begin{eqnarray}\label{roth}
\max_{1 \leq M \leq N} D^*(M, S)\leq D^*(N, P) \leq \max_{1 \leq M \leq N} D^*(M, S)+1.
\end{eqnarray} 
These inequalities will be used to deduce results for two-dimensional point sets from results for one-dimensional sequences and vice versa, see Section 2.

The link between numerical integration and uniformly distributed point sets is provided 
by the Koksma-Hlawka inequality, which bounds the integration error 
of a quasi-Monte Carlo rule by means of the discrepancy of the node sets,
\[
\abs{I_s (f)-\frac{1}{N}\sum_{n=0}^{N-1} f(\bsx_n)}\le V(f)D^*(\{\bsx_0,\bsx_1,\ldots,\bsx_{N-1}\})/N,
\]
where $V(f)$ is the variation of $f$ in the sense of Hardy and Krause (see, e.g., \cite{KN74} for further information).  
However, it should be noted that, considering the huge bounds for the discrepancy 
of low-discrepancy sequences in the usual ranges of $N$, this inequality is unsatisfying and not 
really meaningful (see among others \cite{FL09, FL10}), at least for high dimensions $s$.

For overviews of quasi-Monte Carlo methods and their applications, uniform distribution of point sets, 
and their relations, we refer to the monographs \cite{DP10, DT97, KN74, L09, N92, SJ94}.

\medskip

There are two families of low-discrepancy sequences widely used in QMC methods: Halton sequences \cite{H60} 
and their generalizations, and the so-called $(t,s)$-sequences. In this paper we only consider the second 
family along with $(t,m,s)$-nets, their associated finite point sets of cardinality $b^m$, which we now define in detail. The 
concepts of (digital) $(t,m,s)$-nets and (digital) $(t,s)$-sequences provide very efficient methods to construct point
sets with small star discrepancy. These notions go back to ideas of Sobol' \cite{S67}, Faure \cite{F82}, and Niederreiter \cite{N87}, 
and extensive information on this topic is presented by Niederreiter
in~\cite{N92} (see also \cite{N13} for a recent overview). We first give the general definition of a $(t,m,s)$-net.
\begin{defn}
Let $b\ge 2$, $s\ge 1$, and $0\le t\le m$ be integers. Then a point set $P$ consisting of $b^{m}$ points in
$[0,1)^{s}$ forms a $(t,m,s)$-net in base $b$ if every subinterval
$J=\prod_{j=1}^{s}[a_{j}b^{-d_{j}},(a_{j}+1)b^{-d_{j}})$ of $[0,1)^{s}$, with integers $d_{j}\ge 0$ and integers
$0\le a_{j}< b^{d_{j}}$ for $1\le j\le s$ and of volume $b^{t-m}$, contains exactly $b^{t}$ points of $P$.
\end{defn}
Observe that a $(t,m,s)$-net is extremely well distributed if the quality parameter $t$ is small. 

\begin{examp}\label{Ex1}
A very special, though equally prominent example of a $(0,m,s)$-net in base $b$ is the Hammersley net $\cH_{b,m}$ in dimension $s=2$, consisting 
of $b^m$ points of the form
\[\bsx_n=\left(\phi_b (n),\frac{n}{b^m}\right),\ \ 0\le n\le b^m -1,\]
where $\phi_b$ is the radical inverse function, defined as $\phi_b (n)=\sum_{k=0}^{m-1} \frac{n_k}{b^{k+1}}$ 
for an integer $n$ with base $b$ representation $n=n_0+n_1 b +\cdots + n_{m-1} b^{m-1}$.
\end{examp}

Sometimes we do not consider all $(t,m,s)$-nets in their full generality, 
but restrict ourselves to studying a special construction of these, namely digital $(t,m,s)$-nets over a finite field. 
To this end, as usual,  for $b \ge 2$ we set $\ZZ_b:=\ZZ/b\ZZ$ the residue class ring modulo $b$ equipped with addition 
and multiplication modulo $b$. To keep notation simple, we shall 
sometimes associate the elements of $\ZZ_b$ with the set $\{0,1,\ldots,b-1\}$. Of course, if $p$ is prime $\ZZ_p$ is isomorphic to the field with $p$ elements 
and we do not explicitly distinguish between $\ZZ_p$ and this field. 
 
For prime $p$, a digital $(t,m,s)$-net over $\ZZ_{p}$, which is a special type of a
$(t,m,s)$-net in base $p$, is defined as follows (for a more general definition of $(t,m,s)$-nets over commutative rings, see~\cite{N92}).

For the construction of a digital $(t,m,s)$-net choose $s$ $(m\times m$)-matrices $C_{1},\ldots,C_{s}$ over $\ZZ_{p}$ with
the following property. For each choice of nonnegative integers $d_{1},\ldots,d_{s}$ with
$d_{1}+\cdots+d_{s}=m-t$, the system of the
\begin{center}
\begin{tabular}[t]{l}
first $d_{1}$ rows of $C_{1}$ together with the\\
first $d_{2}$ rows of $C_{2}$ together with the\\
\vdots\\
first $d_{s}$ rows of $C_{s}$
\end{tabular}
\end{center}
is linearly independent over $\ZZ_{p}$. For a fixed $n\in\{0,\ldots,p^{m}-1\}$, let $n$ have base $p$
representation $n=n_{0}+n_{1}p+\cdots +n_{m-1}p^{m-1}$. For $j\in\{1,\ldots,s\}$, multiply the matrix $C_{j}$ by
the vector of digits of $n$ in $\ZZ_{p}^{m}$, which gives
\[C_{j}\cdot(n_{0},\ldots,n_{m-1})^\top=:(y_{1}^{(j)}(n),\ldots,y_{m}^{(j)}(n))^\top\in\ZZ_{p}^{m}.\]
Then we set
\[x_{n}^{(j)}:=\sum_{k=1}^{m}\frac{y_{k}^{(j)}(n)}{p^{k}}.\]
Finally, let $\bsx_{n}:=(x_{n}^{(1)},\ldots,x_{n}^{(s)})$. The point set consisting of the points
$\bsx_{0},\bsx_{1},\ldots,\bsx_{p^{m}-1}$ is called a digital $(t,m,s)$-net over $\ZZ_{p}$ with generating
matrices $C_{1},\ldots,C_{s}$.

\begin{rem}\label{remmultfromright}
Let $P$ be a digital $(t,m,s)$-net over $\ZZ_p$ with generating matrices $C_1,\ldots,C_s$, and let $D$ be a nonsingular 
$(m\times m)$-matrix
over $\ZZ_p$. Then the digital net $Q$ over $\ZZ_p$ that is generated by $C_1\cdot D, C_2\cdot D,\ldots, C_s\cdot D$ is, up to the order of 
the points, the same point set as $P$, since multiplication of the generating matrices by $D$ from the right can be interpreted as a 
re-ordering of the indices $n\in\{0,\ldots, p^m -1\}$ of the points. In particular, $D^\ast (P)= D^\ast (Q)$.   
\end{rem}
\begin{examp}\label{Ex2}
If $b$ is prime, then $\cH_{b,m}$ is a digital $(0,m,2)$-net over $\ZZ_b$ with generating matrices 
\[C_{1}=\begin{pmatrix}
        1&0&\hdots&0\\
        0&\ddots&\ddots &\vdots\\
        \vdots&\ddots&\ddots &0\\
        0&\hdots&0&1
        \end{pmatrix},\
  C_{2}=\begin{pmatrix}
        0&\hdots&0&1\\
        \vdots&\rdots&\rdots &0\\
        0&\rdots&\rdots &\vdots\\
        1&0&\hdots&0
        \end{pmatrix}.\]

\end{examp}

The definition
of a $(t,s)$-sequence is based on $(t,m,s)$-nets and is given in the following.
\begin{defn}
Let $b \ge 2$, $s \ge 1$, and $t \ge 0$ be integers. A sequence $(\bsx_n)_{n \ge 0}$ in $[0,1)^s$ is a
$(t,s)$-sequence in base $b$ if for all $l\ge 0$ and $m > t$ the point set consisting of the points $\bsx_{l
b^m},\ldots, \bsx_{(l+1)b^m-1}$ is a $(t,m,s)$-net in base $b$.
\end{defn}
Again, a $(t,s)$-sequence is particularly well distributed if the quality parameter $t$ is small.
\begin{examp}\label{Ex3}
The probably best known example of a $(0,s)$-sequence in base $b$ is the van der Corput sequence $S_b^\id$ in dimension $s=1$, 
with points of the form
\[S_b^\id(n)=\left(\phi_b (n)\right),\ \ n\ge 0,\]
where the radical inverse function $\phi_b$ is defined as above (see Example \ref{Ex1}).
\end{examp}
A digital $(t,s)$-sequence over $\ZZ_{p}$ ($p$ prime), which is a special type of a $(t,s)$-sequence in base $p$
(cf.~\cite{N92}), is constructed as follows. Choose $s$ $\infty\times\infty$-matrices $C_{1},\ldots,C_{s}$
over $\ZZ_{p}$ such that for any $m>t$ the left upper $m\times m$-submatrices of $C_{1},\ldots,C_{s}$ generate
a digital $(t,m,s)$-net over $\ZZ_{p}$. For $n\ge 0$, let $n$ have base $p$ representation
$n=n_{0}+n_{1}p+\cdots$. For $j\in\{1,\ldots,s\}$, multiply the matrix $C_{j}$ by the vector of digits of $n$,
which yields
\[C_{j}\cdot(n_{0},n_{1},\ldots)^\top=: (y_{1}^{(j)}(n),y_{2}^{(j)}(n),\ldots)^\top,\]
and set
\[x_{n}^{(j)}:=\sum_{k=1}^{\infty}\frac{y_{k}^{(j)}(n)}{p^{k}}.\]
Then the sequence consisting of the points $\bsx_{0},\bsx_{1},\ldots$ with
$\bsx_{n}:=(x_{n}^{(1)},\ldots,x_{n}^{(s)})$ is called a digital $(t,s)$-sequence over $\ZZ_{p}$ and
$C_{1},\ldots,C_{s}$ are its generating matrices. 

A technical requirement on a digital $(t,s)$-sequence is that,
for each $n\ge 0$ and $1\le j\le s$, we have $y_{k}^{(j)}(n)<p-1$ for infinitely many $k$ (cf.~\cite{N92}).

Later, in order to include new important constructions, Niederreiter and Xing \cite{NX96a,NX96} and Tezuka 
\cite{T93} introduced a new definition of $(t,s)$-sequences in arbitrary base $b \ge 2$ using the so-called truncation operator:
\begin{defn}\label{Trunc}
Let $x=\sum_{i=1}^\infty x_i b^{-i}$ be a $b$-adic
expansion of $x\in [0,1]$, where it is allowed that $x_i=b-1$ for all but finitely
many $i$. For every integer $m\geq 1$, the $m$-truncation of $x$ is defined by 
$[x]_{b,m}=\sum_{i=1}^m x_i b^{-i}$ (depending on $x$ via its expansion).
For $\bsx \in I^s $, the notation $[\bsx]_{b,m}$ means that $m$-truncation is applied to each coordinate of $\bsx$.
\end{defn}
\begin{defn}\label{(t,s)broad}
An $s$-dimensional sequence $(\bsx_n)_{n \ge 0}$, with prescribed $b$-adic expansions for each coordinate, 
is a \emph{$(t,s)$-sequence in base $b$ (in the broad sense)} if the subset
$\{[\bsx_n]_{b,m} ; lb^m\le n<(l+1)b^m\}$ is a $(t,m,s)$-net in base $b$
for all integers $l\ge0$ and $m>t.$
\end{defn}
The former $(t,s)$-sequences are now called $(t,s)$-sequences {\it in the narrow sense} 
and the latter (in the sense of Definition \ref{(t,s)broad}) simply $(t,s)$-sequences or sometimes $(t,s)$-sequences 
{\it in the broad sense} (cf. Niederreiter and Xing \cite[Definition 2 and Remark 1]{NX96}). 
These new definitions will be used in the following for one-dimensional sequences.
\begin{examp}
If $b$ is prime, then $S_b^\id$ is a digital $(0,1)$-sequence over $\ZZ_b$ with generating matrix
\[C_{1}=\begin{pmatrix}
        1&0&0&0\\
        0&1&0 &\ddots\\
        0&0&\ddots&\ddots\\
        0&\ddots&\ddots&\ddots
        \end{pmatrix}.\]
\end{examp}

\medskip

The notion of point sets that are {\em digitally permuted} offers a generalization of the notions of nets and sequences. 
These are constructed by a basic variation in the point generating procedures outlined above. To be more precise, 
consider first the one-dimensional case.

Let $\Sy_b$ be the set of all permutations of  $\ZZ_b$. Choose a sequence $\Sigma=(\sigma_k)_{k\ge 0}$ of 
permutations $\sigma_k \in \Sy_b$ and define the sequence 
$S_b^\Sigma$, namely the {\em generalized van der Corput sequence associated with} $\Sigma$ \cite{F81}, by 
\[
S_b^\Sigma(n):=\sum_{k=0}^{\infty} \frac{\sigma_k(n_k)}{b^{k+1}},
\] 
for $n\ge 0$ with base $b$ representation $n=n_0+n_1 b +\cdots$.

Notice that the sequences $S_b^\Sigma$ are $(0,1)$-sequences in the broad sense  \cite[Proposition 3.1]{F07}, the truncation 
from Definition \ref{Trunc} being required to prove this property when $\sigma_k(0)=b-1$ for all sufficiently large $k$.

If, for all $k \ge 0$, $\sigma_k=\sigma$ is constant then we write $S_b^\Sigma=S_b^\sigma$. And if  for all $k \ge 0$ 
we set $\sigma_k=\id$ (the identity), we recover the van der Corput sequence $S_b^\id$ in Example \ref{Ex3} above.

Now, in the $s$-dimensional case, choose $s$ sequences of permutations $\Sigma^{(1)},\ldots, \Sigma^{(s)}$ and 
from an existing point set with points $\bsx_n=(x_n^{(1)},\ldots, x_n^{(s)})$, define the {\em digitally permuted point set} with points 
$\widetilde{\bsx}_n =
(\widetilde{x}_n^{(1)},\ldots, \widetilde{x}_n^{(s)})$, as
\[\widetilde{x}_{n}^{(j)}:=\sum_{k=1}^{\infty}\frac{\sigma_k^{(j)} (x_{n,k}^{(j)})}{b^{k}},\]
where the $x_{n,k}^{(j)}$ are the base $b$ digits of $x_n^{(j)}$.

If the permutations $\sigma$ are of the form $\sigma(x)=f x+ g \pmod b$ with $f \not =0$, we speak of a {\em linear digit scrambling} 
of the point set, and if $f=1$, we speak of a {\em digital shift} $g$ of the point set. Linear digit scramblings are widely 
used in QMC methods to improve the distribution of point sets derived from Halton and $(t,s)$-sequences and hence can
improve approximate computations of integrals (see among others \cite{FL09, FL10}).

\medskip

Regarding the star discrepancy of $(t,m,s)$-nets and $(t,s)$-sequences, it is known by general results of Niederreiter (\cite{N87}) that for any $(t,m,s)$-net $P$ in base $b$
\begin{equation}\label{eqdiscnets}
D^*(P)\le b^t  C(s,b) (\log b^m)^{s-1} +\mathcal{O}\left((\log b^m)^{s-2}\right),
\end{equation} 
with $C(s,b)$ and the implied constant in the $\mathcal{O}$-notation independent of $b^m$. Furthermore, for the star discrepancy of the first $N\ge b^t$ points of a $(t,s)$-sequence $S$ in base $b$,
\begin{equation}\label{eqdiscseq}
D^*(N, S)\le b^t  D(s,b)(\log N)^s +\mathcal{O}\left((\log N)^{s-1}\right),
\end{equation} 
with $D(s,b)$ and the implied constant in the $\mathcal{O}$-notation independent of $N$. Values of $C(s,b)$ and $D(s,b)$ for 
which \eqref{eqdiscnets} and \eqref{eqdiscseq} hold were explicitly given in \cite{N87}, 
and later improved on in \cite{FK13, FL12, K06}.

Concerning lower discrepancy bounds, a famous theorem of Schmidt  \cite{S72} further improved 
with respect to the constant by B\'{e}jian \cite{B82} states that for any sequence $S$ in $[0,1)$ and infinitely many $N$,
\[
D(N, S)\ge 0.12 \log N,
\]
hence achieving the exact order in $N$ for the discrepancy of sequences in one dimension. A similar 
result of B\'{e}jian for the star discrepancy, which states that $D^*(N, S)\ge 0.06 \log N$, was recently improved by Larcher \cite{L13}. 
Furthermore, from \cite[Example 2.2]{KN74}, it follows that
for any two-dimensional point set $P$ of $N$ points, $D^*(N,P)\ge 0.03 \log N$ (see also \cite[Corollary 2.2]{KN74}).

For $s\ge 3$, another famous theorem due to Roth \cite{R54} was recently improved in \cite{BLV08} to 
\[D^*(N, P)\ge c(s) (\log N)^{\frac{s-1}{2}+\delta_s},\] where $c(s)$ is a constant only depending on $s$ and 
not on the point set or its cardinality, and where $\delta_s$ is an unknown constant in $[0,1/2)$.

\medskip

In this paper, we review recent discrepancy estimates for low-dimensional $(t,m,s)$-nets and $(t,s)$-sequences 
(i.e., for $s=1$ and $s=2$). The reason why low-dimensional 
examples of nets and sequences 
have gained much interest by researchers in, particularly, the field of uniform distribution theory, is that special instances (as for 
example, the Hammersley net or the van der Corput sequence) have nice mathematical properties, 
and that a sound understanding of low-dimensional point sets helps in dealing 
with those in arbitrary dimension. 

The rest of the paper is structured as follows. There are three main sections; one (Section \ref{secupperseq}) on 
upper discrepancy bounds for low-dimensional $(t,s)$- and related sequences, one (Section \ref{secuppernets}) on upper discrepancy bounds for 
low-dimensional $(t,m,s)$-nets, and one (Section \ref{seclower}) on lower discrepancy bounds for both 
low-dimensional infinite sequences and low-dimensional finite point sets. 

Section \ref{secupperseq} is divided into several parts. After an introductory part, a few preliminary remarks, and a section on 
general upper bounds for one-dimensional sequences (Section \ref{secgeneral}), Sections \ref{secexact} and \ref{secsbsigma} focus 
on special classes of one-dimensional sequences and contain very precise results on the discrepancy of these. We conclude 
Section \ref{secupperseq} by some remarks on $(t,1)$-sequences (cf. Section \ref{(t,1)-seq}) and two-dimensional sequences (cf. Section 
\ref{2dimseq}).

Section \ref{secuppernets} is organized as follows. In an introductory section we review selected earlier results, and then discuss three 
recent streams of research; one of them is using Walsh functions for deriving discrepancy bounds for certain two-dimensional nets (cf. Section 
\ref{WalshNets}), another one is to use counting arguments for deriving upper discrepancy bounds (cf. Section \ref{seccounting}), 
and one is to relate suitably chosen one-dimensional sequences to generalized Hammersley nets (cf. Section \ref{sechammersley}).

Finally, Section \ref{seclower} contains one section reviewing lower discrepancy bounds for nets (cf. Section \ref{seclowernets}), a
section where we derive a new result (cf. Section \ref{secnew}), and a section discussing lower discrepancy bounds for $(t,s)$-sequences (cf. 
Section \ref{seclowerseq}).

We conclude the paper in Section \ref{secsummary}, where we also state some open questions. 

\section{Upper discrepancy bounds for low-dimensional sequences}\label{secupperseq}

\subsection{Introductory remarks}\label{Sec2Intro}
We leave aside the family of $(n \alpha)$ sequences (see, e.g., \cite{DT97,KN74}), the other great family of one-dimensional low discrepancy sequences,
and focus on $(0,1)$-sequences with a short insight into $(t,1)$-sequences (Section \ref{(t,1)-seq}), 
in relation to $(t,2)$-sequences (Section \ref{2dimseq}).

The family of $(0,1)$-sequences (in the broad sense) contains two large sub-families as shown in 
\cite[Proposition 3.1]{F07}: The family of $S_b^\Sigma$ sequences introduced in 
\cite{F81} and the family of digital $(0,1)$-sequences in (prime, for simplicity) 
base $b$, denoted by $X_b^C$ in the following. In this case, we assume that the generating matrix $C$ has 
the property that for any integer $m \ge 1$ every left upper $m \times m$-submatrix is nonsingular. 
Note that such $X_b^C$ sequences can require the truncation operator (see Definition \ref{Trunc}) since we may have
the digits $y_k(n)=b-1$ for all but finitely many $k$ in the definition. 
An important special case is the case of nonsingular upper triangular (NUT) matrices $C$ 
for which the summation over $k$ in the definition is finite, so that these (so-called) {\em NUT digital $(0,1)$-sequences} do not need the truncation.

Quite recently, Faure and Pillichshammer \cite{FP13} introduced a mixed construction containing both families above.
Such sequences are denoted $X_b^{\Sigma,C}$ and called {\em NUT $(0,1)$-sequences} over 
$\mathbb{Z}_b$ where $b \ge 2$ is an arbitrary base. They are obtained by putting arbitrary permutations from the sequence 
$\Sigma$ in place of the diagonal entries of the NUT matrix $C$. More precisely:
\begin{defn}\label{FaPi13}
For any integer $b\ge 2$, let $\Sigma=(\sigma_k)_{k \ge 0}$ be a sequence of permutations $\sigma_k \in \Sy_b$ and let 
$C=(c_r^k)_{r \ge 0, k \ge r+1}$ be a strict upper triangular matrix with entries in $\ZZ_b$  
(i.e., an upper triangular matrix with naught diagonal entries).
Then, for all integers $n\ge 0$, the $n$-th element of the sequence $X_b^{\Sigma,C}$ is defined by
\begin{equation*}\label{GNUT}
X_b^{\Sigma,C}(n)=\sum_{r=1}^\infty \frac{x^{(r)}(n)}{b^{r}} \quad
{\rm where} \quad x^{(r)}(n)=\sigma_r(n_r)+\sum_{k=r+1}^\infty c_r^k n_k \pmod{b},
\end{equation*} 
in which $n$ has base $b$ representation with digits $n_k$  as in Example \ref{Ex1}.
\end{defn}
If all entries above the diagonal of $C$ are zero, we recover $S_b^\Sigma$ sequences. If the permutations 
$\sigma_k$ are linear digit scramblings with shift $g=0$, we recover classical NUT digital 
$(0,1)$-sequences (with arbitrary base $b$). Sequences $S_b^\Sigma$ and  sequences 
$X_b^C$ are extensively studied in \cite{F81, F05a} but, since their generalization leads formally to the same formulas, 
we will only present these formulas for the case of sequences $X_b^{\Sigma,C}$.

\subsection{Prerequisites}\label{Prereq}
 We first introduce two more notions of discrepancy in one dimension:
\[
D^+(N,X) = \sup_{0\leq \alpha \leq 1} \Delta([0,\alpha),N,X) \mbox{ and }
D^-(N,X) = \sup_{0\leq \alpha \leq 1}(- \Delta([0,\alpha),N,X)).
\]
 The discrepancies $D^+$ and $D^-$ are linked to $D$ and $D^*$ by the relations
\[
D(N,X) = D^+(N,X)+D^-(N,X) \mbox{ and }
D^*(N,X) = \max(D^+(N,X),D^-(N,X)).
\]

Then we need to define the so-called  $\varphi$-functions first introduced in the study of van der Corput sequences \cite{F81}. 
For any $\sigma \in \Sy_b$ (the set of all permutations of $\ZZ_b$),  set 
$$\mathcal{Z}_b^{\sigma}=(\sigma(0)/b,\sigma(1)/b,\ldots ,\sigma(b-1)/b).$$ For $h \in \{0,1,\ldots ,b-1\}$ and $x \in [(k-1)/b,k/b)$ where 
$k \in \{1,\ldots ,b\}$, define 
$$\varphi_{b,h}^{\sigma}(x)=\left\{
\begin{array}{ll}
A([0,h/b);k;\mathcal{Z}_b^{\sigma})-h x & \mbox{ if } 0 \le h \le \sigma(k-1),\\(b-h)x-A([h/b,1);k;\mathcal{Z}_b^{\sigma}) 
& \mbox{ if }\sigma(k-1)< h < b.
\end{array}\right.
$$
Further, the functions $\varphi_{b,h}^{\sigma}$ are extended to the reals by periodicity. Based on $\varphi_{b,h}^{\sigma}$ we now define
 \[
 \psi_b^{\sigma,+}=\max_{0 \le h < b} \varphi_{b,h}^{\sigma},\;\;\;\psi_b^{\sigma,-}=\max_{0 \le h < b} (-\varphi_{b,h}^{\sigma})  
\;\;\mbox{ and } \;\; \psi_b^{\sigma}=\psi_b^{\sigma,+}+\psi_b^{\sigma,-},
 \]
 which appear in the formulas for the discrepancies $D^+$, $D^-$ and $D$.
 
 Moreover, we need further definitions to deal with  (digital or not) NUT $(0,1)$-sequences. The symbol $\uplus$ is used to denote the translation 
(or shift) of a given permutation $\sigma \in \Sy_b$ by an element $t
\in {\mathbb Z}_b$ in the following sense: 
$
(\sigma \uplus t)(i):=\sigma (i)+t \pmod{b}  \mbox{ for \ all } i \in {\mathbb Z}_b,
$
and for any integer $r\ge0$ we introduce the quantity 
\[
\theta_r(N):=\displaystyle\sum_{k=r+1}^\infty c_r^ka_k(N)\pmod{b},
\]
where the $c_r^k$'s are the entries of the matrix $C$ and the $a_k(N)$'s are the digits of $N-1$ in its $b$-adic expansion. Note that
$a_k(N)=0$ for all $k\geq n$ if $1\leq N\leq b^n$, thus $\theta_r(N)=0$ for all $r\geq n-1$ in this case. This quantity determines the 
translated permutations that appear in the formulas for $D^+, D^-$ in Theorem \ref{MThm}.

Lastly, for any $r \ge0$ we define the permutation $\delta_r$ by $\delta_r(i):=c_r^ri \pmod{b}$ for all $i \in \ZZ_b$ and the swapping permutation 
$\tau(i)=b-1-i=(b-1)i+b-1 \pmod{b}$ for all $i \in \ZZ_b$. The name of the permutation $\tau$ comes from the fact that 
$\tau$ swaps the functions $ \psi_b^{\sigma,+}$ and $ \psi_b^{\sigma,-}$ and hence is useful to minimize 
$D^*=\max(D^+,D^-)$ in the asymptotic behavior of the star discrepancy (see below).

\subsection{A general upper bound and two counter-examples}\label{secgeneral}
It is known \cite{F81,F05} that the original van der Corput sequences are the worst distributed 
with respect to the star discrepancy, to the $L_2$ discrepancy and to the diaphony among all 
$S_b^\Sigma$ sequences and among all NUT digital sequences $X_b^C$. This is also true for all 
$(0,1)$-sequences (in the broad sense, see Definition \ref{(t,s)broad}) according to Theorem \ref{thmfk} below. However and 
surprisingly, this property is no longer true for the discrepancy $D$ among all $(0,1)$-sequences 
(in the broad sense) according to our two counter-examples below.
\begin{thm}[Faure, Kritzer] \label{thmfk} The original van der Corput sequences 
in arbitrary base $b\geq 2$ are the worst distributed with respect to the star discrepancy 
among all $(0,1)$-sequences $X_b$ (in the broad sense), that is,
\[
D^*(N,X_b)\leq D^*(N,S_b^\id)=D(N,S_b^\id).
\]
\end{thm}
\begin{rem} {\rm The main idea of the proof was first used by Dick and Kritzer \cite{DK06} 
in the context of two-dimensional Hammersley  point sets. Then Kritzer \cite{K05} proved 
Theorem \ref{thmfk} for $(0,1)$-sequences in the narrow sense using the result for Hammersley 
point sets. Finally, Faure \cite[Theorem 5.1]{F07} proved the theorem  ``in the broad sense"  
using the sequence $S_b^\tau$, whose functions $\psi_b^{\tau,+}$ and $\psi_b^{\tau,-}$ are 
exchanged with functions $\psi_b^{\id,-}$ and $\psi_b^{\id,+}$ associated with $S_b^\id$. 
The good control of discrepancy by means of $\psi$-functions allows a shorter proof.
In the broad sense, we can say that there are two worst sequences with respect to
$D^*$, namely $S_b^\id$ and $S_b^\tau$, while in the narrow sense there is only one, namely $S_b^\id$,
since $S_b^\tau$ is not a sequence in the narrow sense.}
\end{rem}
Now, we give the first counter-example where we show that the original 
van der Corput sequence in base 2 is not the worst distributed sequence with 
respect to $D$ among digital  $(0,1)$-sequences in base $2$ \cite[Theorem 5.2]{F07}.
\begin{thm}[Faure]\label{thm2} Let $C_0$ be the generating matrix in base 2 for which 
all entries are zero except on the diagonal and in the first column (where they are equal to $1$). Then, for any $N \ge 1$,
\[
2D(N,S_2^\id)-\frac{5}{2} \leq D(N,X_2^{C_0}) \leq 2D(N,S_2^\id) \mbox{ and }
\]
\[
D^*(N,S_2^\id)-\frac{3}{2} \leq D^*(N,X_2^{C_0}) \leq D^*(N,S_2^\id)=D(N,S_2^\id).
\]
Moreover the sequence $X_2^{C_0}$ is the worst distributed among all
$(0,1)$-sequences in base $2$ (in the broad sense) with respect to $D$.
\end{thm}
For base $b\ge 3$, we have not found a digital $(0,1)$-sequence with discrepancy $D$ twice greater than that of 
$S_b^\id$. The study of the sequence $X_b^{C_0}$ is, in this case, more complicated and does not 
give the factor 2 like for $b=2$. However, we have a simple construction inspired by the proof of Theorem \ref{thm2} 
which gives the same result \cite[Theorem 5.3]{F07}.
\begin{thm}[Faure]\label{thm3} Let $b\geq 3$ be an integer. Let us define the sequence
$X_b^{\id\tau}=(x_n)_{n\geq 0}$ by
$x_{bk}=S_b^\id(bk), \; x_{bk+1}=S_b^\tau(bk)$ and 
$x_{bk+l}=S_b^\id(bk+l-1)$ if  $2 \leq l \leq b-1$, for all $k\geq 0.$
Then, the sequence $X_b^{\id\tau}$ is a $(0,1)$-sequence (not digital and not in the 
narrow sense), for which
$$2D(N,S_b^\id)-2(b-1) \leq D(N,X_b^{\id\tau}) \leq 2D(N,S_b^\id) \quad{\it and}$$
$$D^*(N,S_b^\id)-(b-1) \leq D^*(N,X_b^{\id\tau}) \leq D^*(N,S_b^\id)=D(N,S_b^\id).$$
Moreover the sequence $X_b^{\id\tau}$ is the worst distributed among all
$(0,1)$-sequences in base $b$ (in the broad sense) with respect to $D$.
\end{thm}

\subsection{Exact formulas for the discrepancies  of NUT $(0,1)$-sequences}\label{secexact}
We now turn to results on NUT $(0,1)$-sequences $X_b^{\Sigma,C}$ as defined in Section \ref{Sec2Intro}, Definition \ref{FaPi13}. 
Here, we restrict ourselves to the discrepancies $D^+$, $D^-$ (and so $D^*$) and $D$. 
Similar formulas exist for the $L_2$ discrepancy and the diaphony (see \cite{FP13}).
\begin{thm}[Faure, Pillichshammer] \label{MThm}
With the notation introduced in Sections \ref{Sec2Intro} and \ref{Prereq}, we have, for all integers $b \ge 2$ and  $N \ge 1$,
\begin{equation*}\label{D+}
D^+(N,X_b^{\Sigma,C})=\sum_{j=1}^\infty
\psi_b^{\sigma_{j-1}\uplus\theta_{j-1}(N),+} \left ( \frac{N}{b^j} \right ),
\end{equation*}
\begin{equation*}\label{D-}
D^-(N,X_b^{\Sigma,C})=\sum_{j=1}^\infty \psi_b^{\sigma_{j-1}\uplus
\theta_{j-1}(N),-}\left ( \frac{N}{ b^j} \right ),
\end{equation*}
\begin{equation*}\label{D}
D(N,X_b^{\Sigma,C})=\sum_{j=1}^\infty \psi_b^{\sigma_{j-1}}
\left ( \frac{N}{b^j} \right ).
\end{equation*}
\end{thm}

For sequences $X_b^{\Sigma,C}$, the formula for $D$ depends only on the permutations $\sigma_r$ 
and not on the  entries (above the diagonal) of $C$.  For sequences $X_b^{C}$ in \cite{F05}  
the formula for $D$ depends only on the permutations $\delta_r$ (associated with the diagonal entries of $C$) 
and not on the  entries above the diagonal of $C$. This remarkable feature shows that NUT (0,1)-sequences 
having the same sequence of permutations $\Sigma$  (or $\Delta=(\delta_r)_{r \ge 0}$ for $X_b^{C}$ sequences) have 
the same extreme discrepancy $D$. In this case, the studies on the asymptotic behavior of $D$ for generalized van der Corput
sequences  $S_b^\Sigma$ (see \cite{F81, F92, F08}) apply,  
especially to NUT digital (0,1)-sequences.  The same remark is also valid for the diaphony, see \cite{F05, FP13}.

On the other hand, the formulas  for $D^+$, $D^-$, and hence for $D^*$, 
involve the quantity $\theta_{j-1}(N)$ which depends on $N$ via its $b$-adic
expansion and on the generating NUT matrix $C$ via its entries above the
diagonal; this dependence is a big handicap for the precise study of the asymptotic
behavior of the star discrepancy of NUT (0,1)-sequences. 
As far as we know, the only result  available is that of Pillichshammer for $X_2^C$ \cite{P04} 
where $C$ is the matrix for which all entries are equal to $1$. Indeed, Pillichshammer showed the following result, which 
can be 
deduced from the 
study of digital $(0,m,2)$-nets in base 2 by Larcher and Pillichshammer \cite{LP03} (see Section \ref{WalshNets}).
\begin{thm}[Pillichshammer]
 For the star discrepancy of the sequence $X_2^C$ described above it is true that
$$ 0.2885\ldots =\frac{1}{5\log 2}\le \limsup_{N\To\infty} \frac{D^\ast (N,X_2^C)}{\log N}\le \frac{5099}{22528\log 2}=0.3265\ldots\ .$$
\end{thm}

\subsection{Upper bounds for one-dimensional low discrepancy sequences}


\paragraph{Reminders on the asymptotic behavior of  $S_b^\Sigma$ sequences.}\label{secsbsigma}
We only recall here two main theorems going back to 1981 useful for our short 
review of results on upper bounds for one-dimensional low discrepancy sequences \cite[Th\'eor\`emes 2 and 3]{F81}.

For an infinite sequence $X$ in $[0,1)$, set 
\[ \rho(X):=\limsup_{N\rightarrow\infty}\frac{D(N,X)}{\log N} \;\mbox{ and }\;
\rho^*(X):=\limsup_{N\rightarrow\infty}\frac{D^*(N,X)}{\log N}.
\]
\begin{thm}[Faure]\label{thm2F81} 
Let $\sigma \in \Sy_b$ and let $\Sigma=(\sigma)$ be constant (so that $S_b^\Sigma=S_b^\sigma$). Then
\[ \rho(S_b^\sigma)=\frac{\alpha_b^{\sigma}}{2 \log b} \quad
\mbox{ where } \quad
\alpha_b^{\sigma}=\inf_{n\geq 1} \frac{1}{n}\sup_{x\in [0,1]} \sum_{j=1}^n
\psi_b^{\sigma}\left (\frac{x}{b^j}\right ) \quad \mbox{and}
\]
\[ D(N,S_b^\sigma)\leq \displaystyle\frac{\alpha_b^\sigma }{\log b} \log N+\alpha_b^\sigma+2 \mbox{ for all } N \ge 1.
\]
\end{thm}
The previous result concerns the extreme discrepancy. We now turn to the star discrepancy, 
where best results are obtained with special sequences of permutations $\Sigma$ using the swapping permutation $\tau$ 
introduced at the end of Section \ref{Prereq} and recalled below:.
\begin{thm}[Faure]\label{thm3F81} Let $\tau \in \Sy_b$ be the permutation defined by $\tau(k)=b-k-1$ for all $k \in \ZZ_b$ and let 
$\mathcal{A}$ be the subset of $\NN_0$ defined by  
$\mathcal{A}=\bigcup_{H=1}^\infty \mathcal{A}_H$ with $\mathcal{A}_H=\{H(H-1),\ldots,H^2-1\}$.
For any permutation $\sigma \in \Sy_b$, let $\overline{\sigma}:=\tau \circ \sigma$ and let
$$\Sigma^\sigma_\mathcal{A}=(\sigma_r)_{r \ge 0} := (\sigma,\overline{\sigma},\sigma,\sigma,\overline{\sigma},
\overline{\sigma},\sigma,\sigma,\sigma,\overline{\sigma},\overline{\sigma},\overline{\sigma},\ldots)$$
be the sequence of permutations defined by $\sigma_r=\sigma$ if $r \in \mathcal{A}$ and $\sigma_r=\overline{\sigma}$ if $r \notin \mathcal{A}$. Then
\begin{equation}\label{f_jnt}
\rho^*(S_b^{\Sigma^\sigma_\mathcal{A}})=\frac{\alpha_b^{\sigma,+}+\alpha_b^{\sigma,-}}{2 \log b},
\end{equation}
where
$$\alpha_b^{\sigma,+}=\inf_{n\geq 1} \frac{1}{n}\sup_{x\in [0,1]} \sum_{j=1}^n
\psi_b^{\sigma,+}\left (\frac{x}{b^j}\right ) \quad \mbox{ and } \quad  
\alpha_b^{\sigma,-}=\inf_{n\geq 1} \frac{1}{n} \sup_{x\in [0,1]} \sum_{j=1}^n
\psi_b^{\sigma,-}\left (\frac{x}{b^j}\right ).
$$
\end{thm}


\paragraph{Updated review of results on one-dimensional low discrepancy sequences.}
For reasonably small $b$, the constants $\alpha_b^{\sigma,+}$ and $\alpha_b^{\sigma,-}$ are not difficult to compute and for the identity 
$\id$, in which case $\psi_b^{\id,-}=0$, it is even possible to find them explicitly. We have 
\[
\frac{\alpha_b^{\id}}{\log b}=\frac{b-1}{4\log b} \mbox{ if } b \mbox{ is odd,\quad and \quad} \frac{\alpha_b^{\id}}{\log b}= \frac{b^2}{4(b+1)\log b} 
\mbox{ if } b \mbox{ is even}.
\]
These constants are the worst possible leading constants for the discrepancies $D$ and $D^*$ of $S_b^\Sigma$ sequences since for any sequence of 
permutations $\Sigma$, we have \cite[Section 5.5.4 ]{F81}
\[
D(N,S_b^\Sigma)\leq D(N,S_b^{\id})=D^*(N,S_b^{\id}) \mbox{ for all } N\ge1,\]
(i.e., the original van der Corput sequence $S_b^{\id}$ is the worst distributed among the $S_b^\Sigma$ sequences with respect to the discrepancy and the 
star discrepancy).

Concerning the sequences in Theorem \ref{thm3F81}, we obtain
\[
\rho(S_b^{\Sigma^\id_\mathcal{A}})=\rho^*(S_b^{\Sigma^\id_\mathcal{A}})=\frac{\alpha_b^{I,+}}{2 \log b}=\frac{b-1}{8 \log b} \; 
\mbox{ if } b \mbox{ is odd, and}\]
\[
\rho(S_b^{\Sigma^\id_\mathcal{A}})=\rho^*(S_b^{\Sigma^\id_\mathcal{A}})=\frac{\alpha_b^{I,+}}{2 \log b}=\frac{b^2}{8 (b+1)\log b}\; 
\mbox{ if } b \mbox{ is even}.
\]
In these formulas, base $b=2$ is interesting: in this case the swapping permutation $\tau$ reads as $\tau(k)=k+1 \pmod{2}$ and so 
$\tau$ is the digital shift\!$\pmod{2}$. Hence in 2007, Kritzer, Larcher and Pillichshammer \cite[Section 5]{KLP07} re-discovered Theorem 
\ref{thm3F81} in the special case of base $2$ and the corresponding constant $1/(6 \log 2)$ above. 

Notice that in general we have $\alpha_b^\sigma \le \alpha_b^{\sigma,+}+\alpha_b^{\sigma,-}$. Hence we cannot infer any relation between 
a general lower bound on $\rho^*(S_b^{\Sigma^\sigma_\mathcal{S}})$, where $\mathcal{S}$ is an arbitrary subset of $\NN_0$, and the upper bound on 
$\rho^*(S_b^{\Sigma^\sigma_\mathcal{A}})$ from \eqref{f_jnt}. Only if for the permutation $\sigma$ we have either $\psi_b^{\sigma,+}=0$ or 
$\psi_b^{\sigma,-}=0$, in which case we get $\alpha_b^\sigma =\alpha_b^{\sigma,+}+\alpha_b^{\sigma,-}$, we obtain that 
$\alpha_b^\sigma / (2 \log b)$ is the best possible lower bound for the star discrepancy of sequences $S_b^\Sigma$ with 
$\Sigma \in \{\sigma,\overline{\sigma}\}^{\NN_0}$. In other words, we have 
\[
\inf_{\Sigma \in \{\sigma,\tau\circ\sigma\}^{\NN_0}} \rho^*(S_b^\Sigma)= \frac{\alpha_b^\sigma}{2 \log b}
\]
for any $\sigma \in \Sy_b$ such that $D^*(S_b^\sigma)=D(S_b^\sigma)$. Recall that for any $\Sigma = (\sigma_r)_{r \ge 0}$ we have 
$D^*(S_b^\Sigma)=D(S_b^\Sigma)$ if and only if  $\psi_b^{{\sigma_r},+}=0$ for all $r \ge 0$ or $\psi_b^{{\sigma_r},-}=0$  for all $r \ge 0$. 
See \cite[Corollaire 2, p. 160]{F81}) and \cite[Section 5.2]{FP13} for more details..

The smallest extreme discrepancies currently known, a new record,  
were recently obtained by Ostromoukhov~\cite{O09} in bases 84 and 60, after a lot 
of computations using the method worked out in \cite[Section 5]{F81}: there exist  
permutations $\sigma_0  \in \Sy_{84}$  and  $\sigma_1 \in \Sy_{60}$ such that
\[
\rho(S_{84}^{\Sigma^{\sigma_0}_\mathcal{A}})=\frac{130}{83 \log 84}=0.3534\ldots \mbox{ and }
\rho^*(S_{60}^{\Sigma^{\sigma_1}_\mathcal{A}})=\frac{32209}{35400 \log 60}=0.2222\ldots,
\]
improving preceding results in \cite[Theorem 1.2]{F92} and \cite[Th\'eor\`eme 5]{F81} with permutations 
$\sigma_2 \in \Sy_{36}$ and $\sigma_3 \in \Sy_{12}$ that give $\rho(S_{36}^{\sigma_2})=0.3667\ldots$ and 
$\rho^*(S_{12}^{\Sigma^{\sigma_3}_\mathcal{A}})=0.2235\ldots$. It is interesting to note that $\psi_{12}^{\sigma_1,-} \neq0$ whereas 
$\psi_{60}^{\sigma_0,-}=0$. This last property is quite remarkable with regard to the multitude of permutations involved in the computational search for 
$(60,\sigma_0)$ among all pairs $(b,\sigma)$.


\subsection{A general upper bound for $(t,1)$-sequences}\label{(t,1)-seq}
As announced at the beginning of Section \ref{Sec2Intro}, we now turn to a slight 
generalization of Theorem \ref{thmfk}, see \cite[Section 2]{FL12}.
\begin{thm} [Faure, Lemieux] \label{thmFL} For any base $b$, the original van der 
Corput sequences are the worst distributed with respect to the star discrepancy among all $(t,1)$-sequences
$X_b^t$ (in the broad sense), i.e., for all $N \ge 1$,
\[
D^*(b^tN,X_b^t)\leq b^tD^*(N,S_b^\id)=b^tD(N,S_b^\id).
\]
\end{thm}
As a corollary, we obtain that for any $(t,1)$-sequence $X_b^t$ (in the broad sense) and for any integer $N \ge 1$,
\begin{align*}
D^*(b^t N,X_b^t)\leq& b^t \left (\frac{b-1}{4\log b} \log N + \frac{b-1}{4}+2 \right )\mbox{ if } b \mbox{ is odd, and }\\
D^*(b^t N,X_b^t)\leq& b^t \left (\frac{b^2}{4(b+1)\log b} \log N + \frac{b^2}{4(b+1)}+2 \right ) \mbox{ if } b \mbox{ is even}. 
\end{align*}
These upper bounds in one dimension complete the upper bounds in dimension $s \ge 2$ recently obtained by Faure and Kritzer, see Section \ref{2dimseq} below.


\subsection{Two-dimensional sequences}\label{2dimseq} Regarding two-dimensional sequences, the currently best known general upper discrepancy 
bound was first shown for digital sequences in base $2$ in \cite{P03}, and then for arbitrary $(t,2)$-sequences in \cite{DK06}:
\begin{thm}[Dick, Kritzer]\label{thmup02seq}
 For the star discrepancy of the first $N$ points of an arbitrary $(t,2)$-sequence $X_b^t$ in base $b$, it is true that
 \[N D_N^\ast (X_b^t)\le \begin{cases} \frac{b^t}{16}\frac{b^2 (b-1)^2}{(b^2-1)(\log b)^2} (\log N)^2 + \mathcal{O}(\log N)  
                     &\mbox{if $b$ is even,}\\ \\
                     \frac{b^t}{16}\frac{(b-1)^2}{(\log b)^2} (\log N)^2 + \mathcal{O} (\log N)
                     &\mbox{if $b$ is odd,} \end{cases} \]
where the implied constants in the $\mathcal{O}$-notation do not depend on $N$.
\end{thm}
\begin{rem}{\rm
Theorem \ref{thmup02seq} has just been generalized to arbitrary dimensions $s \ge 2$ by Faure and Kritzer \cite[Theorem 2 and Corollary 2]{FK13}. 
The proofs follow the approach consisting first in doing the study of 
$(t,m,s)$-nets for which formulas are proved with the help of a recursion based on preliminary counting lemmas; 
then going from $(t,m,s)$-nets to $(t,s)$-sequences results from a classical relation already used by Sobol', 
Faure and Niederreiter (see \cite[Lemma 4.1]{N87}). There is another approach dealing directly with 
$(t,s)$-sequences and using an adaptation of Atanassov's method for Halton sequences \cite{FL12}, but the 
bounds obtained in this way are greater by a factor of about 2 (even though the constants  hidden in the 
$\mathcal{O}$-notation are smaller in this approach). 
Another interesting remark is that  Corollary 2 in \cite{FK13} and Theorem \ref{thmup02seq} are also true in one dimension thanks to Theorem \ref{thmFL}.
}
\end{rem}


\section{Upper discrepancy bounds for low-dimensional nets}\label{secuppernets}
Let us now turn to upper discrepancy bounds for low-dimensional finite point sets. 


\subsection{Introductory remarks}\label{Sec3Intro}
For the case $s=1$, there exists a closed formula for the discrepancy of a finite 
point set in $[0,1)$, so one usually  does not explicitly study the discrepancy of $(t,m,1)$-nets. 
Indeed, assume that $x_0,\ldots,x_{N-1} \in [0,1)$ such that, without 
loss of generality, $x_0\le x_1\le \cdots \le x_{N-1}$. Then
\begin{equation}\label{eqexplicit}
D^{\ast} (\{x_0,\ldots,x_{N-1}\})=\frac{1}{2} + N \max_{0\le n\le N-1}\abs{x_n - \frac{2n+1}{2N}}.
\end{equation}
Due to this formula, the derivation of discrepancy bounds for finite one-dimensional point sets is obsolete. However,  if one wants to study infinite point sets 
$S$ in $[0,1)$ and derive bounds on $D^\ast(N,S)$ that hold for all or at least infinitely many $N$, this is still a non-trivial problem, see Section \ref{secupperseq}.

Upper discrepancy bounds for $(t,m,s)$-nets were, in a very general form, shown in 
\cite{N87, N92}, where the formula \eqref{eqdiscnets} was proved. In particular, 
Niederreiter also presents discrepancy bounds for dimension $s=2$ in \cite{N87,N92}, namely
\begin{thm}[Niederreiter]\label{thmniedspec} 
The star discrepancy of a $(t,m,2)$-net $P$ in base $b$ satisfies
  \begin{equation}\label{eqnied} 
  D^{\ast}(P) \le \left\lfloor \frac{b-1}{2} (m-t) + \frac{3}{2}\right\rfloor b^t.
  \end{equation}
\end{thm}

Regarding further, special discrepancy bounds for $(t,m,2)$-nets, there exist results that can be grouped into 
two larger streams of research. On the one hand, approaches by Faure used explicit formulas for the local discrepancy function $\Delta$, 
see \cite[Lemma 3]{F08a}), which make it possible to deduce discrepancy bounds for low-dimensional point sets from the study 
of low-dimensional sequences. On the other hand, 
initiated by Larcher and Pillichshammer in 2003, cf.~\cite{LP03}, a series of papers dealt with tightening discrepancy bounds 
for low-dimensional nets by using Walsh series.


\subsection{Using Walsh functions to obtain bounds for $(0,m,2)$-nets in base $2$}\label{WalshNets}
Regarding the results of Larcher and Pillichshammer in~\cite{LP03}, 
these are relevant for parts of the remainder of this paper (see Theorem \ref{thmnew} in Section 3.2), so let us shortly recall them here. 
In~\cite{LP03}, the discrepancy of digital $(0,m,2)$-nets over $\ZZ_2$ is analyzed by the means of Walsh functions. 
Indeed, let $P$ be a digital $(0,m,2)$-net over $\ZZ_2$. Without loss of generality (see Remark~\ref{remmultfromright}), we assume that
$P$ is generated by matrices
\[C_{1}=\begin{pmatrix}
        1&0&\hdots&0\\
        0&\ddots&\ddots &\vdots\\
        \vdots&\ddots&\ddots &0\\
        0&\hdots&0&1
        \end{pmatrix},\
  C_{2}=\begin{pmatrix}
        c_{1,1}&c_{1,2}&\hdots&c_{1,m}\\
        c_{2,1}&c_{2,2}&\hdots&c_{2,m}\\
        \vdots&\vdots&\vdots &\vdots\\
        c_{m,1}&c_{m,2}&\hdots&c_{m,m}
        \end{pmatrix}=
        \begin{pmatrix} \bsc_{1}\\ \bsc_{2}\\ \vdots\\ \bsc_{m}\end{pmatrix}.\]
If we would like to estimate $D^\ast (P)$, then it is often sufficient to 
restrict oneself to studying the local discrepancy function $\Delta (J,P)$ only for certain well-chosen intervals $J$. In \cite{LP03},  
the authors considered the local discrepancy function $\Delta$ of a digital net, evaluated at
$(\eta,\beta)\in [0,1)^2$ for $m$-bit numbers $\eta$ and $\beta$, i.e., numbers with base 2 representation
\[\eta=\frac{\eta_1}{2}+\frac{\eta_2}{2^2}+\cdots + \frac{\eta_m}{2^m},\ 
  \ \ \beta=\frac{\beta_1}{2}+\frac{\beta_2}{2^2}+\cdots + \frac{\beta_m}{2^m}.\]
As pointed out in~\cite{LP03},
\[\abs{D^\ast (P)-2^m \max_{\eta,\beta\  m-\mathrm{bit}}\abs{\Delta([0,\eta)\times [0,\beta),P)}}\le 2-1/2^m,\]
which means that one can in many cases resort oneself to studying the discrepancy of $P$ only for $m$-bit numbers, 
without significant deviations from the precise value of the discrepancy. 

We introduce some further notation. We write $\bseta:=(\eta_1,\ldots,\eta_m)^\top$ and $\bsbeta:=(\beta_1,\ldots,\beta_m)^\top$. Let
\[\bsgamma=\bsgamma (\eta,\beta)=(\gamma_{1}, \gamma_{2}, \ldots, \gamma_{m})^{\top}=C_2\cdot\bseta+\bsbeta,\]
and let $\bsgamma (u)$ denote the vector consisting of the first $u$ components of $\bsgamma$.

Furthermore, we write, for $1\le u\le m$, 
  \[C_2 ' (u)= \begin{pmatrix}
        c_{1,m-u+1}&c_{2,m-u+1}&\hdots&c_{u,m-u+1}\\
        c_{1,m-u+2}&c_{2,m-u+2}&\hdots&c_{u,m-u+1}\\
        \vdots&\vdots&\vdots &\vdots\\
        c_{1,m}&c_{2,m}&\hdots&c_{u,m}
        \end{pmatrix}^{-1},\]
which exists according to the $(0,m,2)$-net property of $P$. 

In the following, we denote by $(\cdot | \cdot)$ the usual inner product.
Moreover, we denote by $\norm{x}$ the distance to the nearest integer of a real number $x$. Notice that 
$\norm{x}=\psi_2^\id(x)$, the $\psi$-function in base 2 associated with identity in  Section \ref{Prereq}.

In \cite{LP03}, the following formula for the local discrepancy $\Delta(\eta,\beta)$ was shown,
\begin{equation}\label{eqformuladelta}
\Delta(\eta,\beta)=\sum_{u=0}^{m-1}\norm{2^u \beta} (-1)^{(\bsc_{u+1} \vert \bseta)}(-1)^{(\bsgamma (u)
  \vert C_2' (u)(c_{u+1,m-u+1},\ldots,c_{u+1,m})^\top)} \frac{(-1)^{\eta_{m-u}}- (-1)^{\eta_{m+1-j(u)}}}{2},
\end{equation}
where
\[m(u)=\begin{cases}
         0 \ \  \mbox{if}\  u=0,&\\
         0 \ \  \mbox{if}\ (\bsgamma (u)\vert C_2' (u)\bse_1)=1,&\\
         \max\{1 \le j\le u:\ (\bsgamma (u)\vert C_2' (u) \bse_i)=0,\  i=1,\ldots,j\}& \mbox{otherwise},
       \end{cases}
\]
and $j(u)=u-m(u)$. Here $\bse_i$ denotes the $i$-th unit vector in $\ZZ_2^u$. Furthermore, we set $\eta_{m+1}:=0$ and
$(\bsgamma (u)\vert C_2' (u)(c_{u+1,m-u+1},\ldots,c_{u+1,m})^\top):=0$ for $u=0$. 

The formula \eqref{eqformuladelta} is a powerful tool. Indeed, by means of \eqref{eqformuladelta}, 
it was shown in \cite{LP03} that the following improvement of Theorem~\ref{thmniedspec} holds. 
\begin{thm}[Larcher, Pillichshammer]\label{thmlp}
 The star discrepancy of a digital $(0,m,2)$-net $P$ over $\ZZ_2$ satisfies
\begin{equation}\label{eqlarpil}
D^{\ast}(P) \le m/3 + 19/9.
\end{equation}
\end{thm}
This result was further sharpened, and extended to arbitrary two-dimensional nets, in the following Theorem \ref{thmdk} in \cite{DK06}.


\subsection{Using a counting argument to obtain bounds for $(t,m,2)$-nets}\label{seccounting}
 Surprisingly, the proof of Theorem \ref{thmdk} does not need any of the technical tools mentioned above, but is based on a counting argument. 
\begin{thm}[Dick, Kritzer]\label{thmdk}
The star discrepancy of an arbitrary $(t,m,2)$-net $P$ in base $b$ satisfies
\begin{equation}\label{eqdk}
D^{\ast}(P) \le b^t D^{\ast} (\cH_{b,m-t}) + b^{t}.
\end{equation} 
\end{thm}
The bound in Theorem \ref{thmdk} is effectively useful, as $D^{\ast} (\cH_{b,m})$ can be computed exactly by using formulas provided 
by DeClerck in \cite{DC86}. In particular, this implies that \eqref{eqdk} is indeed an improvement of \eqref{eqlarpil} for the case $t=0, b=2$, 
and also an improvement of \eqref{eqnied} for arbitrary choices of $t$ and $b$. 
Furthermore, the bound in \eqref{eqdk} can be seen as an extension of the corollary to Theorem \ref{th1H} below and is, up to the constant $b^{t}$, sharp. 

If we focus on the case $t=0$, we know by Equation~\eqref{eqdk} that the Hammersley point set is basically the $(0,m,2)$-net in base $b$ with 
the highest star discrepancy. This observation sparked interest in the question whether there are examples of other $(0,m,2)$-nets with a 
significantly lower star discrepancy than $\cH_{b,m}$. This question can be answered positively, and, not surprisingly, 
digitally permuted Hammersley nets play a crucial role in this context, since they are connected to 
digitally permuted $(0,1)$-sequences (see Section \ref{sechammersley} below).


\subsection{Using results on $S_b^\Sigma$ sequences to obtain bounds for Hammersley nets}\label{sechammersley}
We first define what we mean by digitally permuted Hammersley nets associated with $S_b^\Sigma$ sequences. In order to match the traditional 
definition of arbitrary (shifted or not) Hammersley point sets whose points are ``$m$-bits'', we restrict the infinite sequence of permutations 
$\Sigma$ to permutations such that $\sigma_r(0)=0$ for all $r\ge m$, for instance $\bssigma:=(\sigma_0,\ldots,\sigma_{m-1},\id ,\id, \id, \ldots)$. 
Then the {\it generalized two-dimensional Hammersley point set in base $b$} consisting of $b^m$ points associated with $\bssigma$ is defined by
$$\cH_{b,m}^{\bssigma}:=\left \{ \left (S_b^{\bssigma}(n), \frac{n}{b^m} \right );0\leq n \leq b^m-1 \right \}.$$
Notice that the behavior of $\cH_{b,m}^{\bssigma}$ only depends on the finite sequence $(\sigma_0,\ldots,\sigma_{m-1})$ which we identify with $\bssigma$ from now on (see \cite[Definition 2]{F08a} for more details).
If we choose in the above definition $\sigma_j=\id$ for all $j$, then we obtain the classical two-dimensional Hammersley point 
set in base $b$, $\cH_{b,m}^{\bsid}=\cH_{b,m}$.

The main results concerning the star discrepancy of two-dimensional Hammersley point 
sets are of two kinds: some give exact formulas including complementary terms (\cite{DC86,HZ69,LP03}) 
and the others give formulas for the leading terms within an error not computable, usually lower than a small 
additive constant (\cite{DK06, F86, F08, K06, KLP07}). For the sake of brevity, we will only refer to these latter 
results in the generalizations we are going to give in the following. Theorems \ref{th1H} to \ref{th5H} below stem from \cite{F08a}.

First, we make the tight link between $\cH_{b,m}^{\bssigma}$ and $S_b^{\bssigma}$ via $\psi$-functions more precise :
\begin{thm}[Faure]\label{th1H}
For any integer $m\ge 1$ and any $\bssigma=(\sigma_0,\ldots,\sigma_{m-1})$ we have, with some $c_m \in [0,2]$,
\[
D^*(\cH_{b,m}^{\bssigma})= \max \left (\max_{1 \le n \le b^m}\sum_{j=1}^m \psi_b^{\sigma_{j-1},+} \bigg( \frac{n}{b^j} \bigg)\ , 
\ \max_{1 \le n \le b^m} \sum_{j=1}^m \psi_b^{\sigma_{j-1},-} \bigg( \frac{n}{b^j} \bigg)\right)+c_m.
\]
\end{thm}
This formula, valid for an arbitrary base $b$, is the analog of \cite[Lemma 1]{KLP07} in base $2$. It permits all extensions 
to bases $b$ of results on base $2$ belonging to the second kind mentioned above.

As a first corollary, we obtain that $\cH_{b,m}^{\bsid}$ is the worst Hammersley net among generalized ones within a constant less than 2:
For any integer $m\ge 1$ and any $\bssigma=(\sigma_0,\ldots,\sigma_{m-1})$ we have, with some $c_m \in [-2,2]$,
\begin{equation}\label{eqf1}
D^*(\cH_{b,m}^{\bssigma}) \le D^*(\cH_{b,m}^{\bsid})+c_m.
\end{equation}

\paragraph{Swapping with the identical permutation.}

First let us consider the following sequence $\bsi \bstau=
(\overbrace{\id,\ldots,\id}^{\frac{m}{2}}, \overbrace{\tau, \ldots \tau}^{\frac{m}{2}})$ if $m$ is even and $\bsi \bstau=(\overbrace{\id,\ldots,\id}^{\frac{m-1}{2}}, 
\overbrace{\tau, \ldots \tau}^{\frac{m+1}{2}})$, if $m$ is odd, like Kritzer did in base 2 \cite{K06a}. Applying Theorem \ref{th1H}, we 
can easily extend his result \cite[Theorem 3.2 and Proposition 3.1]{K06a} to arbitrary bases:

\begin{thm}[Faure]\label{th2H}
For any integer $m\ge 1$ we have, with some $c_m \in [0,3]$,
$$\mbox{if } b \mbox{ is odd:} \quad D^*(\cH_{b,m}^{\bsi \bstau})=\left\{
\begin{array}{ll}
\displaystyle \frac{b-1}{8}m+c_m & \mbox{ if } m \mbox{ is even}\\
\displaystyle \frac{b-1}{8}(m+1)+c_m & \mbox{ if } m \mbox{ is odd,}
\end{array}\right.
$$
$$\mbox{if } b \mbox{ is even:} \quad D^*(\cH_{b,m}^{\bsi \bstau})=\left\{
\begin{array}{ll}
\displaystyle \frac{b^2}{8(b+1)}m+c_m & \mbox{ if } m \mbox{ is even}\\
\displaystyle \frac{b^2}{8(b+1)}(m+1)+c_m & \mbox{ if } m \mbox{ is odd.}
\end{array}\right.
$$
\end{thm}
The interval for $c_m$ could be reduced.
Of course, we recover the result of Kritzer \cite{K06a} and Kritzer, Larcher and Pillichshammer \cite{KLP07}, 
with the same sequence $\bsi \bstau$, in the case of $b=2$. 
The best constant is obtained for $b=3$ with $1/(4\log 3)=0.227\ldots$, whereas for $b=2$ we only have $1/(6\log 2)=0.240\ldots$.

Now, we show that the choice of the sequence $\bsi \bstau$ in the set $\{\id,\tau \}^m$ of sequences $\bssigma=(\sigma_0, \ldots, \sigma_{m-1})$ 
is best possible in the sense that the leading terms in Theorem \ref{th2H} cannot be made smaller whatever 
the $\sigma_{j-1} \in \{\id,\tau \}$, $1 \le j \le m$, are. 
\begin{thm}[Faure]\label{th3H}
For any integer $m\ge 1$ and any $\bssigma \in \{\id,\tau \}^m$ we have
$$\lim_{m \rightarrow \infty} \frac{D^*(\cH_{b,m}^{\bssigma})}{\log b^m}\ge \frac{b-1}{8\log b} \quad \mbox{ if } b \mbox{ is odd and }$$
$$\lim_{m \rightarrow \infty} \frac{D^*(\cH_{b,m}^{\bssigma})}{\log b^m}\ge \frac{b^2}{8(b+1)\log b}\quad \mbox{ if } b \mbox{ is even.}$$
\end{thm}
In base $2$, Theorem \ref{th3H} has been shown in \cite{KLP07} by other arguments involving more computations. 
Another question raised and solved in base 2 in \cite{KLP07} is the following: 
\vskip 3mm 
\noindent
{\it ``Is the star discrepancy $D^*(\cH_{b,m}^{\bssigma})$ independent of the distribution of $\id$ and $\tau$ in 
the sequence $\bssigma=(\sigma_0,\ldots,\sigma_{m-1}) \in \{\id,\tau \}^m$ and does only depend on the number of $\id$ and $\tau$?"}
\vskip 3mm
In arbitrary base, the answer is {\it No} like in base 2, with the same counter-example as in \cite{KLP07}, the sequence $(\id, \tau, \id, \tau, \ldots, \id, \tau)$.
\begin{thm}[Faure]\label{th4H}
For any even integer $m\ge 2$, let $\widetilde{\bsi \bstau}=(\id, \tau, \id, \tau, \ldots,\id, \tau) \in \{\id,\tau \}^m$.
Then, with some $c_m \in [0,3]$, we have
$$D^*(\cH_{b,m}^{\widetilde{\bsi \bstau}})=\left\{
\begin{array}{ll}
\displaystyle\frac{(b-1)(b+2)}{8(b+1)}m+c_m & \mbox{if } b \mbox{ is odd,}\\
\displaystyle\frac{b^3}{8(b^2+1)}m+c_m & \mbox{if } b \mbox{ is even}.
\end{array}\right.$$
These constants are greater than $(b-1)/8$ and $b^2/(8(b+1))$, hence the answer to the question above is {\it No}.
\end{thm}
For $b=2$ we recover the result of \cite[end of Section 4]{KLP07} with the constant $1/5$. This result has been known for long since Halton and Zaremba 
\cite{HZ69} obtained, in 1969, exact formulas for $D^*(\cH_{2,m}^{\bsid})$ and $D^*(\cH_{2,m}^{\widetilde{\bsi \bstau}})$ after a lot of technical computations (see \cite[Sections 1 and 4]{LP03} for comments).
\medskip

\paragraph{Swapping with an arbitrary permutation.} In this section, we fix an arbitrary permutation $\sigma$ of $\{0,1,\ldots ,b-1\}$ and consider sequences produced by swapping $\sigma$ with 
$\tau$, like we did in Theorem \ref{th2H}, to obtain sequences $\bssigma=(\sigma_0,\ldots,\sigma_{m-1}) \in \{\sigma,\overline{\sigma} \}^m$. 
The situation is not so clear as with the identity and we will only consider sequences $\bssigma \overline{\bssigma}=
(\overbrace{\sigma,\ldots,\sigma}^{\frac{m}{2}}, \overbrace{\overline{\sigma}, \ldots \overline{\sigma}}^{\frac{m}{2}})$ if $m$ is even and 
$\bssigma \overline{\bssigma}=(\overbrace{\sigma,\ldots,\sigma}^{\frac{m-1}{2}}, \overbrace{\overline{\sigma}, \ldots \overline{\sigma}}^{\frac{m+1}{2}})$ if 
$m$ is odd. This choice permits to improve on the discrepancy, but until now we have not been able to prove it is the best, like with the identity.
\begin{thm}[Faure]\label{th5H} For any integer $m \ge 1$ we have, with some $c_m \in [-1,4]$,
\[
D^*(\cH_{b,m}^{\bssigma\overline{\bssigma}})=\left \{
\begin{array}{ll}
\displaystyle\frac{\alpha_b^{\sigma,+}+\alpha_b^{\sigma,-}}{2}m+c_m & \mbox{ if } m \mbox{ is even}\\
\displaystyle\frac{\alpha_b^{\sigma,+}+\alpha_b^{\sigma,-}}{2}(m+1)+c_m & \mbox{ if } m \mbox{ is odd.}
\end{array}\right.
\]
\end{thm}
Of course, we recover the constants of Theorem \ref{th2H} when $\sigma=\id$. Obtained as a by-product of Theorem \ref{thm3F81}, the best result coming from 
\cite[Th\'eor\`eme 5]{F81}, with $b=12$, $\sigma_3$ and  constant $0.2235\ldots$, was recently improved by Ostromoukhov \cite{O09} with $b=60$, 
$\sigma_1$ and  constant $0.2222\ldots$, a bit better than $b=3$ and $\sigma=\id$ with $0.227\ldots$ (see the comments following Theorem \ref{thm3}). 
Even if such improvements seem small, they concern the leading constants in discrepancy formulas and we think it is more important to improve 
on these constants rather than searching for exact formulas or to reduce the complementary terms $c_m$ in estimations. 

\section{Lower discrepancy bounds for low-dimensional point sets}\label{seclower}

\subsection{Lower discrepancy bounds for nets}\label{seclowernets}

We now survey lower bounds on the star discrepancy of low-dimensional nets. As in Section \ref{secuppernets}, we focus on the case of dimension $s=2$, 
as the case $s=1$ is essentially answered by Equation (\ref{eqexplicit}). 

Generally speaking, it is a lot more challenging to provide tight lower discrepancy bounds than upper bounds.
The only general result that holds for all $(t,m,2)$-nets is the aforementioned bound of B\'ejian \cite{B82} (see the end of Section \ref{Intro}), 
from which it results that 
any finite point set $P$ of $N$ points in $[0,1)^2$ satisfies 
\begin{equation}
D^\ast(N,P) \ge 0.03 \log N.
\end{equation}
In view of the bounds \eqref{eqlarpil}, \eqref{eqdk}, \eqref{eqf1}, many researchers have tried to sharpen lower bounds for at least some 
well-chosen subclasses of $(t,m,2)$-nets, and there has been a major focus on $(0,m,2)$-nets. 

For example, it was shown in \cite{KLP07} that a digitally shifted Hammersley net $\cH_{2,m}$ in base 2 always has a discrepancy such that 
$D^\ast(\cH_{2,m}) \ge m/6 +c$, where $c$ is a constant independent of $m$. 
Faure \cite{F08} showed a similar result in this vein, proving that 
it is true for any digitally permuted Hammersley net $\cH_{b,m}^{\bssigma}$ in base $b$ that 
\[D^\ast(\cH_{b,m}^{\bssigma})\ge \left(3/4 - \sqrt{3b-1}/(2b)\right)m + c,\]
where $c$ is some positive constant that does not depend on $m$.

\subsection{A new result on the discrepancy of digital $(0,m,2)$-nets}\label{secnew}

In this section, we show a new lower discrepancy bound that holds for certain digital $(0,m,2)$-nets over $\ZZ_2$. 
For $m\in\NN$, we write $m_0:= \left\lfloor \frac{m}{2} \right\rfloor$. Within this section, we consider digital $(0,m,2)$-nets over $\ZZ_2$ 
generated by 
\begin{equation}\label{eqspecialchoice}
C_1=I_{m,m},\ \ \  C_2 =\begin{pmatrix} A &\vline & B\\
                       \hline
                       C &\vline & D  \end{pmatrix},
\end{equation}
where $I_{m,m}$ is the $m\times m$ identity matrix, and
where $A$ a nonsingular $m_0\times m_0$-matrix over $\ZZ_2$. 
We show the following new result.
\begin{thm}\label{thmnew}
 Let $P$ be a digital $(0,m,2)$-net generated by two generating matrices as in \eqref{eqspecialchoice}. Then it is true that
\[D^\ast (P)\ge m/12 + c,\]
where $c$ is a constant independent of $m$. 
\end{thm}
\begin{rem}{\rm
 We remark that, even though the bound in Theorem \ref{thmnew} is weaker than the lower bound for digitally shifted Hammersley nets from \cite{KLP07} 
 mentioned before, the result in Theorem \ref{thmnew} covers a relatively large class of digital $(0,m,2)$-nets.}
\end{rem}
\begin{proof}
We use the approach of Larcher and Pillichshammer \cite{LP03} that was 
summarized in Section \ref{secuppernets} and the same notation as there. We make a specific choice of two base 2 $m$-bit numbers $\eta$ and $\beta$. Let
\[\beta^{(0)}=(b_{m_0},\underbrace{0,0,\ldots\ldots\ldots\ldots,0}_{m-m_0\ \mathrm{components}}),\]
where 
\[b_{m_0}=\begin{cases}
       \underbrace{1,0,1,0,\ldots,1,0,1,0}_{m_0\ \mathrm{components}}& \mbox{if}\ m_0\ \mbox{is even},\\
       \underbrace{1,0,1,0,\ldots 1,0,1}_{m_0\ \mathrm{components}}& \mbox{if}\ m_0\ \mbox{is odd}.\\
      \end{cases}
\]

For the following, denote by $\bszero_{k,l}$ the $k\times l$ zero matrix. We define an $m$-bit number $\eta^{(0)}$ such that
\[\begin{pmatrix}
   A & \vline & B\\
   \bszero_{m-m_0,m_0}&\vline & I_{m-m_0, m-m_0}
  \end{pmatrix} \cdot \bseta^{(0)} = \begin{pmatrix} \delta_1\\ \delta_2 \\ \vdots \\ \delta_m \end{pmatrix}=\bsdelta,
\]
where the first and the last $m_0$ components, respectively, of $\bsdelta$ satisfy
\[\delta_{m-u}\oplus 1=\delta_{u+1} = \beta_{u+1}^{(0)}\ \mbox{for}\ 0\le u\le m_0-1,\]
and where $\oplus$ denotes addition modulo 2.
Note that $\eta^{(0)}$ can certainly be chosen in such a way due to the assumptions made on the matrix $A$ in Equation \eqref{eqspecialchoice}.

Using the notation in Section \ref{secuppernets}, note that the choice of $\eta^{(0)}$ and $\beta^{(0)}$ guarantees that
$\bsgamma (u)=\bszero$
for $1\le u\le m_0-1$, which implies that $m(u)=u$ and $j(u)=0$ for  $0\le u\le m_0 -1$. Note furthermore, that $\norm{2^u \beta^{(0)}}=0$ for
all $u\ge m_0$. Hence the discrepancy function $\Delta$ evaluated at $(\eta^{(0)},\beta^{(0)})$ simplifies to
\[\Delta(\eta^{(0)},\beta^{(0)})=\sum_{u=0}^{m_0-1} \norm{2^u \beta^{(0)}} (-1)^{\beta^{(0)}_{u+1}} \frac{(-1)^{\beta^{(0)}_{u+1}+ 1}-1}{2}.\]
However,
\[(-1)^{\beta^{(0)}_{u+1}} \frac{(-1)^{\beta^{(0)}_{u+1}+ 1}-1}{2}=\begin{cases} 
  -1 &\mbox{if}\ \beta^{(0)}_{u+1}=0,\\
  0 &\mbox{if}\ \beta^{(0)}_{u+1}=1.
  \end{cases}\]
Thus,
\[\abs{\Delta(\eta^{(0)},\beta^{(0)})}=\sum_{\substack{u=1\\u\ \mathrm{odd}}}^{m_0-1} \norm{2^u \beta^{(0)}}.\]
Note that, for the case of $m_0$ being even,
\begin{eqnarray}\label{eqm0even}
\sum_{\substack{u=1\\u\ \mathrm{odd}}}^{m_0-1} \norm{2^u \beta^{(0)}}=
  \sum_{\substack{u=1\\u\ \mathrm{odd}}}^{m_0-1}\norm{ 2^u\sum_{k=1}^{m_0/2}\frac{1}{2^{2k-1}}}
&=&\sum_{\substack{u=1\\u\ \mathrm{odd}}}^{m_0-1}\quad\sum_{k=1}^{(m_0-1-u)/2}\frac{1}{2^{2k}}\nonumber\\
&=&\sum_{\substack{u=1\\u\ \mathrm{odd}}}^{m_0-1}\frac{1}{3}\left(1-\frac{1}{2^{m_0-1-u}}\right)\nonumber\\
&=&\sum_{k=1}^{m_0/2}\frac{1}{3}\left(1-\frac{1}{2^{m_0-1-(2k-1)}}\right)\nonumber\\
&=&\frac{m_0}{6} +\frac{4}{9}\left(\frac{1}{2^{m_0}}-1\right),
\end{eqnarray}
and, for the case of $m_0$ being odd, we can show in a similar way as in the derivation of \eqref{eqm0even},
\begin{equation}\label{eqm0odd}
\sum_{\substack{u=1\\u\ \mathrm{odd}}}^{m_0-1} \norm{2^u \beta^{(0)}}=\frac{m_0}{6} +\frac{1}{9}\left(\frac{1}{2^{m_0}}-1\right).
\end{equation}
Now note that $m_0$ is of order $m/2$, and hence we conclude from \eqref{eqm0even} and \eqref{eqm0odd} that
\[\abs{\Delta(\eta^{(0)},\beta^{(0)})}\ge m/12 + c,\]
where $c$ is some constant independent of $m$. The result follows.
\end{proof}

\subsection{Lower discrepancy bounds for sequences}\label{seclowerseq}
\paragraph{Lower discrepancy bounds for one-dimensional sequences.} 
In this part, we give an application of Theorem~\ref{MThm} and show best possible lower bounds on the star discrepancy of NUT $(0,1)$-sequences. 
This study is motivated by a best possible lower bound on the star discrepancy of digitally shifted van der Corput sequences in base 2 shown 
in \cite{KLP07} and the question whether this bound remains true also for digitally shifted NUT digital sequences in base 2.

Theorems \ref{MThm} and \ref{thm3F81} are two main ingredients of the following result which leads to best possible lower bounds for 
large sub-families of NUT $(0,1)$-sequences to be stated afterwards. 
\begin{thm}[Faure, Pillichshammer] \label{LowBd}
For any integer $b \ge 2$, let $\sigma \in \Sy_b$ and let $C$ be a strict upper 
triangular matrix with entries in $\ZZ_b$. Then, for any subset $\mathcal{S}$ of $\NN_0$, we have
$$
D^*(N,X_b^{\Sigma^\sigma_\mathcal{S},C}) \ge \frac{1}{2}D(N,S_b^\sigma) \;\;\mbox{ and hence}\;\;\; \rho^*(X_b^{\Sigma^\sigma_\mathcal{S},C}) 
\ge \frac{\alpha_b^\sigma}{2 \log b},
$$
where $\Sigma^\sigma_\mathcal{S}=(\sigma_r)_{r \ge 0}$ with $\sigma_r=\sigma$ if $r \in \mathcal{S}$ and $\sigma_r=\tau \circ \sigma$ if 
$r \notin \mathcal{S}$, and $X_b^{\Sigma^\sigma_\mathcal{S},C}$ is a NUT $(0,1)$-sequence.
\end{thm}
We start with  a best possible lower bound for NUT $(0,1)$-sequences $X_b^{\Sigma,C}$ associated with sequences of permutations 
$\Sigma \in \{\sigma,\tau\circ\sigma\}^{\NN_0}$ for which $\sigma$ gives permuted van der Corput sequences with $D=D^*$.  
\begin{cor}\label{BestLowBd}
Let $\mathcal{C}_{\mathrm{SUT}}$ be the set of all strict upper triangular matrices and let $\sigma \in \Sy_b \mbox{ such that }\\ D^*(S_b^\sigma)=D(S_b^\sigma)$. 
Then
$$
\inf_{\substack{\Sigma \in \{\sigma,\tau\circ\sigma\}^{\NN_0} \\ C \in \mathcal{C}_{\mathrm{SUT}}}}
\rho^*(X_b^{\Sigma,C}) = \frac{\alpha_b^\sigma}{2 \log b}.
$$
\end{cor}
Besides the identity $\id$, it is not difficult to find permutations satisfying the condition $D^*(S_b^\sigma)=D(S_b^\sigma)$. 
Further, a systematic computer search performed by F. Pausinger (IST Austria, personal communication) 
has given 26, 58, 340, and 1496 such permutations in bases 6, 7, 8, and 9, respectively.
\medskip

The case of identity in Corollary~\ref{BestLowBd} is of special interest because $\alpha_b^\id$ is explicitly known for any integer $b \ge 2$.

\begin{cor}\label{infbX}
With the notation of Corollary \ref{BestLowBd}, we obtain
$$
\inf_{b \ge 2} \inf_{\substack{\Sigma \in \{\id,\tau\}^{\NN_0} \\ C \in \mathcal{C}_{\mathrm{SUT}}}}
\rho^* (X_b^{\Sigma,C}) = \frac{1}{4 \log 3}=0.2275\ldots.
$$
\end{cor}
This result can be seen as the analog for NUT $(0,1)$-sequences of the best possible 
lower bound for the star discrepancy of $(n\alpha)$ sequences obtained by Dupain and S\'{o}s \cite{DS84}, with $\rho^*((n \sqrt{2}))=0.2836\ldots$ .
We see that NUT $(0,1)$-sequences yield a much smaller value of $\rho^*$.

\medskip
Finally, we consider digitally permuted NUT digital sequences by means of linear digit scramblings. Such sequences, denoted $Z_b^{\Pi,C}$, are defined as follows: let $\mathcal{C}_{NUT}^1$ be the set of NUT matrices $C$ such that all the diagonal entries $c_r^r=1$ and let $\Pi=(\pi_r)_{r \ge 0} \in \Sy_b^{\NN_0}$ be a sequence of linear digit scramblings.
Then, for any $n \ge 0$,
$$
Z_b^{\Pi,C}(n)=\sum_{r=0}^\infty \frac{\pi_r(x_{n,r})}{b^{r+1}} \quad
{\rm with } \quad x_{n,r}=\sum_{k=r}^\infty c_r^k n_k \pmod{b},
$$
where the $n_k$ are the base $b$ digits of $n$. We have the following analog of Corollary \ref{infbX}:
\begin{cor}\label{Picase}
Let $Z_b^{\Pi,C}$ be a linearly digit scrambled NUT digital $(0,1)$-sequence associated with  $C \in \mathcal{C}_{\mathrm{NUT}}^1$ and 
$\Pi=(\pi_r)_{r \ge 0} \in \{\id, \tau\}^{\NN_0}$. Then we have
$$
\inf_{b \ge 2} \inf_{\substack{\Pi \in \{\id,\tau\}^{\NN_0} \\ C \in \mathcal{C}_{\mathrm{NUT}}^1}}
\rho^*(Z_b^{\Pi,C}) = \frac{1}{4 \log 3}=0.2275\ldots.
$$
\end{cor}
\noindent
Notice that  $\id$ and $ \tau$ are the only linear digit scramblings satisfying  $D^*(S_b^\pi)=D(S_b^\pi)$.

\medskip
In the case $b=2$, Corollary \ref{Picase} permits to answer the question evoked at the beginning: ``Is it true that the constant 
$1/(6 \log 2)$ is best possible for any {\em digitally shifted NUT digital sequence in base 2}, as it is the case for any digitally 
shifted van der Corput sequence according to \cite[Corollary~4]{KLP07}?" Taking into account that, in base $2$, $\tau$ is the nonzero 
shift and the diagonal entries of $C$ are all equal to $1$ we can answer this question in the affirmative:

\begin{cor}\label{dsdNUTseq}
We have
$$
\inf_{\substack{\Delta \in \ZZ_2^{\NN_0} \\ C \in \mathcal{C}_{\mathrm{NUT}}}}
\rho^*(Z_2^{\Delta,C})=\frac{1}{6 \log 2}.
$$
\end{cor}
For more information on the context of Theorem \ref{LowBd} and its corollaries we refer to \cite[Section 5]{FP13},
where an overview of this topic is given.

\paragraph{Lower discrepancy bounds for two-dimensional sequences.} 
Regarding lower bounds for $(t,2)$-sequences, only very little is known, except for one example by Faure and Chaix 
in \cite{FC96}, who were able to obtain the exact order of the star discrepancy  
for a $(0,2)$-sequence $S_{\mathrm{Sob}}$ in base 2 first introduced by Sobol' \cite{S67}:
\begin{thm}[Faure, Chaix]\label{LowBdSob}
The digital $(0,2)$-sequence in base 2 generated by the identity matrix and the Pascal matrix$\mod{2}$, denoted $S_{\mathrm{Sob}}$, 
satisfies the inequality
\[
\frac{1}{24(\log 2)^2} \le \limsup_{N \rightarrow \infty} \frac{D^\ast (N,S_{\mathrm{Sob}})}{(\log N)^2}\cdot
\]
\end{thm}
In combination with Theorem \ref{thmup02seq}, this is the only case of a low discrepancy sequence in dimension greater than one 
for which the exact order of discrepancy is known.
Moreover, based on thorough numerical experiments that permitted to find the subsequence leading to their lower bound, Faure and Chaix 
stated the conjecture that the inequality above should actually be an equality, i.e.,
$
\limsup_{N \rightarrow \infty} \frac{D^\ast(N,S_{\mathrm{Sob}})}{(\log N)^2} = \frac{1}{24(\log 2)^2}\cdot
$

\section{Conclusion}\label{secsummary}
In this survey, we have illustrated that there has 
been a considerable history of results on discrepancy bounds for low-dimensional $(t,m,s)$-nets, $(t,s)$-sequences, and related point sets. 
We have summarized 20 theorems, a large part of them dealing with two-dimensional nets, an equally large part dealing with one-dimensional sequences, 
and further results on two-dimensional sequences. All results on one-dimensional sequences in Section \ref{secupperseq} stem from the initial study 
\cite{F81} on generalized van der Corput sequences, Theorem \ref{MThm} being the foremost new generalization for these sequences. Regarding two-dimensional nets in Section 
\ref{secuppernets}, we have discussed several different approaches; two theorems are obtained by counting arguments (Theorems \ref{thmniedspec} and \ref{thmdk}), two others result from Walsh 
function analysis of discrepancy (Theorems \ref{thmlp} and \ref{thmnew}) and the remaining five theorems stem from the study of generalized van der 
Corput sequences (Theorems \ref{th1H}--\ref{th5H}). While these results deal with the precise study of a special class of nets (namely generalized Hammersley nets), 
the previous ones concern arbitrary $(t,m,2)$-nets (Theorems \ref{thmniedspec} and \ref{thmdk}) or digital $(0,m,2)$-nets (Theorem \ref{thmlp} and Theorem \ref{thmnew}, 
which is the only previously unpublished result of this paper). 

Finally, two main results deal with two-dimensional sequences: Theorem \ref{thmup02seq} on upper bounds for arbitrary 
$(t,2)$-sequences has been recently extended to arbitrary $(t,s)$-sequences \cite{FK13}, and Theorem \ref{LowBdSob} on lower bounds remains the 
only exception for which the exact order is attained in dimension $s>1$. 

Most of the results mentioned in this paper have been obtained by methods of number theory and algebra, and the 
precise analysis of the discrepancy of the point sets, even though they are ``only'' one- or two-dimensional, is very challenging. The recent results
on the discrepancy of $(t,m,s)$-nets and $(t,s)$-sequences in arbitrary dimension $s$ in \cite{FK13}, 
which are partly obtained by an inductive argument on the dimension $s$, demonstrate that it may be very crucial to have excellent discrepancy bounds for low-dimensional 
examples, as the better the low-dimensional starting point, the better the results obtained inductively can be expected to be.

We end this paper by stating two selected open problems that would be interesting to be solved in the near future.
\vskip 2pt
$\bullet$ Find other two-dimensional sequences than that in Theorem \ref{LowBdSob} having the ``correct'' order of star discrepancy. Natural candidates are 
$(0,2)$-sequences in arbitrary bases and two-dimensional Halton sequences, for instance that in bases 2 and 3. This open problem seems to be a very challenging task.
\vskip 2pt
$\bullet$ Find an exact formula  for the discrepancy function of one-dimensional digital sequences or two-dimensional digital nets in base $b$, 
i.e., extend Formula (\ref{eqformuladelta}) from \cite{LP03} to other bases $b>2$ . Such a formula is extended to arbitrary  bases for digital 
$(0,1)$-sequences generated by NUT matrices in \cite{F05a}, but until now, no analog exists for generating matrices having nonzero entries below the diagonal. 
In relation with this question, we refer to Theorem \ref{thm2}, where a generating matrix having nonzero entries below the diagonal leads to a surprising result. 
Further investigations on such matrices could help to make progress in the understanding of digital nets and sequences and their distribution properties.

\section*{Acknowledgements}

The authors would like to thank G.~Larcher and F.~Pillichshammer for suggestions and remarks, 
and the referee for valuable comments improving the consistency of this article.

\begin{small}
\noindent\textbf{Authors' addresses:}\\

\medskip

\noindent Henri Faure\\
Institut de Math\'{e}matiques de Luminy, UMR 6206 CNRS\\
163 Av. de Luminy, case 907, 13288 Marseille cedex 9, France\\
E-mail: \texttt{faure@iml.uni-mrs.fr}

\medskip

\noindent Peter Kritzer\\
Department of Financial Mathematics, Johannes Kepler University Linz\\ 
Altenbergerstr.~69, 4040 Linz, Austria\\
E-mail: \texttt{peter.kritzer@jku.at}
\end{small}

\end{document}